\documentclass[11pt]{amsart}
\usepackage{amsmath,amssymb}
\newtheorem{theorem}{Theorem}
\newtheorem{corollary}[theorem]{Corollary}
\newtheorem{lemma}[theorem]{Lemma}
\newtheorem{proposition}[theorem]{Proposition}

\begin{document}

\title[Einstein manifolds with nonnegative isotropic curvature]{Einstein manifolds with nonnegative isotropic curvature are locally symmetric}
\author{Simon Brendle}
\address{Department of Mathematics \\
                 Stanford University \\
                 Stanford, CA 94305}
\thanks{This work was supported by the National Science Foundation under grants DMS-0605223 and DMS-0905628.}

\maketitle

\section{Introduction}

The study of Einstein manifolds has a long history in Riemannian geometry. An important problem, first studied by M.~Berger \cite{Berger1}, \cite{Berger2}, is to classify all Einstein manifolds satisfying a suitable curvature condition. For example, if $(M,g)$ is a compact Einstein manifold of dimension $n$ whose sectional curvatures lie in the interval $(\frac{3n}{7n-4},1]$, then $(M,g)$ has constant sectional curvature (see \cite{Besse}, Section 0.33). A famous theorem of S.~Tachibana \cite{Tachibana} asserts that a compact Einstein manifold with positive curvature operator has constant sectional curvature. Moreover, Tachibana proved that a compact Einstein manifold with nonnegative curvature operator is locally symmetric. M.~Gursky and C.~LeBrun \cite{Gursky-LeBrun} have obtained interesting results on four-dimensional Einstein manifolds with nonnegative sectional curvature. Another result in this direction was established by D.~Yang \cite{Yang}. 

We now describe a curvature condition which was introduced by M.~Micallef and J.D.~Moore \cite{Micallef-Moore}. To that end, let $(M,g)$ be a compact Riemannian manifold of dimension $n \geq 4$. We say that $(M,g)$ has positive isotropic curvature if 
\begin{align*} 
&R(e_1,e_3,e_1,e_3) + R(e_1,e_4,e_1,e_4) \\ 
&+ R(e_2,e_3,e_2,e_3) + R(e_2,e_4,e_2,e_4) \\ 
&- 2 \, R(e_1,e_2,e_3,e_4) > 0 
\end{align*} 
for all orthonormal four-frames $\{e_1,e_2,e_3,e_4\} \subset T_p M$. Moreover, we say that $(M,g)$ has nonnegative isotropic curvature if 
\begin{align*} 
&R(e_1,e_3,e_1,e_3) + R(e_1,e_4,e_1,e_4) \\ 
&+ R(e_2,e_3,e_2,e_3) + R(e_2,e_4,e_2,e_4) \\ 
&- 2 \, R(e_1,e_2,e_3,e_4) \geq 0 
\end{align*} 
for all orthonormal four-frames $\{e_1,e_2,e_3,e_4\} \subset T_p M$. It was shown in \cite{Brendle-Schoen1} that positive isotropic curvature is preserved by the Ricci flow in all dimensions (see also \cite{Nguyen}). This fact plays a central role in the proof of the Differentiable Sphere Theorem (cf. \cite{Brendle-Schoen1}, \cite{Brendle-Schoen2}, \cite{Brendle}).

M.~Micallef and M.~Wang showed that a four-dimensional Einstein manifold with nonnegative isotropic curvature is locally symmetric (see \cite{Micallef-Wang}, Theorem 4.4). In this paper, we extend the results of Micallef and Wang to higher dimensions:

\begin{theorem}
\label{main.theorem}
Let $(M,g)$ be a compact Einstein manifold of dimension $n \geq 4$. If $(M,g)$ has positive isotropic curvature, then $(M,g)$ has constant sectional curvature. Moreover, if $(M,g)$ has nonnegative isotropic curvature, then $(M,g)$ is locally symmetric.
\end{theorem}

We note that H.~Seshadri \cite{Seshadri} has obtained an interesting partial classification of manifolds with nonnegative isotropic curvature. 

We now give an outline of the proof of Theorem \ref{main.theorem}. Let $(M,g)$ be a compact Einstein manifold with nonnegative isotropic curvature. Moreover, suppose that $(M,g)$ is not locally symmetric. After passing to the universal cover if necessary, we may assume that $M$ is simply connected. We now consider the holonomy group of $(M,g)$.

If $\text{\rm Hol}(M,g) = SO(n)$, then $(M,g)$ has positive isotropic curvature. We then show that $(M,g)$ has constant sectional curvature. The proof uses the maximum principle, as well as an algebraic inequality established in \cite{Brendle-Schoen1}.

If $n = 2m \geq	4$ and $\text{\rm Hol}(M,g) = U(m)$, then $(M,g)$ is a K\"ahler-Einstein manifold with positive orthogonal bisectional curvature. It then follows from work of S.~Goldberg and S.~Kobayashi \cite{Goldberg-Kobayashi} that $(M,g)$ is isometric to $\mathbb{CP}^m$ up to scaling. 

If $n = 4m \geq 8$ and $\text{\rm Hol}(M,g) = \text{\rm Sp}(m) \cdot \text{\rm Sp}(1)$, then $(M,g)$ is a quaternionic-K\"ahler manifold. By a theorem of Alekseevskii (cf. \cite{Besse}, Section 14.41), the curvature tensor of $(M,g)$ can be written in the form $R = R_1 + \kappa \, R_0$, where $R_1$ has the algebraic properties of a hyper-K\"ahler curvature tensor, $R_0$ is the curvature tensor of $\mathbb{HP}^m$, and $\kappa$ is a constant. Since $(M,g)$ has nonnegative isotropic curvature, we have $R_1(X,JX,X,JX) < \kappa$ for all points $p \in M$ and all unit vectors $X \in T_p M$. Using the maximum principle, we are able to show that $R_1(X,JX,X,JX) \leq 0$ for all points $p \in M$ and all unit vectors $X \in T_p M$. From this, we deduce that $R_1$ vanishes identically. Consequently, the manifold $(M,g)$ is isometric to $\mathbb{HP}^m$ up to scaling. From this, the assertion follows.

M.~Berger \cite{Berger3} has shown that every quaternionic-K\"ahler manifold with positive sectional curvature is isometric to $\mathbb{HP}^m$ up to scaling. C.~LeBrun and S.~Salamon \cite{LeBrun-Salamon} have conjectured that a quaternionic-K\"ahler manifold $(M,g)$ with positive scalar curvature is necessarily locally symmetric. The results in this paper imply that no counterexample to the LeBrun-Salamon conjecture can have nonnegative isotropic curvature.

Part of this work was carried out during a visit to ETH Z\"urich, Switzerland. I would like to thank Professor Michael Struwe and Professor Tristan Rivi\`ere for inspiring discussions. Finally, I am grateful to the referee for useful comments on an earlier version of this paper.

\section{Preliminary results}

Let $V$ be a finite-dimensional vector space equipped with an inner product. An algebraic curvature tensor on $V$ is a multi-linear form $R: V \times V \times V \times V \to \mathbb{R}$ satisfying 
\[R(X,Y,Z,W) = -R(Y,X,Z,W) = R(Z,W,X,Y)\] 
and 
\[R(X,Y,Z,W) + R(Y,Z,X,W) + R(Z,X,Y,W) = 0\] 
for all vectors $X,Y,Z,W \in V$.

Let $\{e_1, \hdots, e_n\}$ be an orthonormal basis of $V$. Moreover, suppose that $R$ and $S$ are two algebraic curvature tensors on $V$. We define an algebraic curvature tensor $B(R,S)$ on $V$ by 
\begin{align*} 
&B(R,S)(X,Y,Z,W) \\ 
&= \frac{1}{2} \sum_{p,q=1}^n \big [ R(X,Y,e_p,e_q) \, S(Z,W,e_p,e_q) + R(Z,W,e_p,e_q) \, S(X,Y,e_p,e_q) \big ] \\ 
&+ \sum_{p,q=1}^n \big [ R(X,e_p,Z,e_q) \, S(Y,e_p,W,e_q) + R(Y,e_p,W,e_q) \, S(X,e_p,Z,e_q) \big ] \\ 
&- \sum_{p,q=1}^n \big [ R(X,e_p,W,e_q) \, S(Y,e_p,Z,e_q) + R(Y,e_p,Z,e_q) \, S(X,e_p,W,e_q) \big ]
\end{align*} 
for all vectors $X,Y,Z,W \in V$. Finally, for each algebraic curvature tensor $R$, we define $Q(R) = B(R,R)$. 

The following result is purely algebraic: 

\begin{proposition} 
\label{Q.terms}
Let $V$ be a vector space of dimension $n \geq 4$ which is equipped with an inner product. Let $R$ be an algebraic curvature tensor on $V$ with nonnegative isotropic curvature. Finally, suppose that $\{e_1,e_2,e_3,e_4\}$ is an orthonormal four-frame in $V$ satisfying 
\begin{align*} 
&R(e_1,e_3,e_1,e_3) + R(e_1,e_4,e_1,e_4) \\ 
&+ R(e_2,e_3,e_2,e_3) + R(e_2,e_4,e_2,e_4) \\ 
&- 2 \, R(e_1,e_2,e_3,e_4) = 0. 
\end{align*} 
Then 
\begin{align*} 
&Q(R)(e_1,e_3,e_1,e_3) + Q(R)(e_1,e_4,e_1,e_4) \\ 
&+ Q(R)(e_2,e_3,e_2,e_3) + Q(R)(e_2,e_4,e_2,e_4) \\ 
&- 2 \, Q(R)(e_1,e_2,e_3,e_4) \geq 0. 
\end{align*} 
\end{proposition} 

\textbf{Proof.} 
This was shown in \cite{Brendle-Schoen1} (see Corollary 10 in that paper). \\

The term $Q(R)$ arises naturally in the evolution equation for the curvature tensor under Ricci flow (cf. \cite{Hamilton1}, \cite{Hamilton2}). In the special case of Einstein manifolds, we have the following well-known result:

\begin{proposition}
\label{curvature.identity}
Let $(M,g)$ be a Riemannian manifold with $\text{\rm Ric}_g = \rho \, g$. Then the Riemann curvature tensor of $(M,g)$ satisfies 
\[\Delta R + Q(R) = 2\rho \, R.\] 
\end{proposition}

\textbf{Proof.} 
It follows from Lemma 7.2 in \cite{Hamilton1} that 
\begin{align*}
&(\Delta R)(X,Y,Z,W) + Q(R)(X,Y,Z,W) \\ 
&= (D_{X,Z}^2 \text{\rm Ric})(Y,W) - (D_{X,W}^2 \text{\rm Ric})(Y,Z) \\ &- (D_{Y,Z}^2 \text{\rm Ric})(X,W) + (D_{Y,W}^2 \text{\rm Ric})(X,Z) \\ 
&+ \sum_{k=1}^n \text{\rm Ric}(X,e_k) \, R(e_k,Y,Z,W) + \sum_{k=1}^n \text{\rm Ric}(Y,e_k) \, R(X,e_k,Z,W) 
\end{align*}
for all vector fields $X,Y,Z,W$. Since $\text{\rm Ric}_g = \rho \, g$, we conclude that 
\[(\Delta R)(X,Y,Z,W) + Q(R)(X,Y,Z,W) = 2\rho \, R(X,Y,Z,W),\] 
as claimed. \\

Finally, we shall need the following result:

\begin{proposition} 
\label{borderline}
Let $(M,g)$ be a compact Einstein manifold of dimension $n \geq 4$ with nonnegative isotropic curvature. Then the set of all orthonormal four-frames $\{e_1,e_2,e_3,e_4\}$ satisfying 
\begin{align*} 
&R(e_1,e_3,e_1,e_3) + R(e_1,e_4,e_1,e_4) \\ 
&+ R(e_2,e_3,e_2,e_3) + R(e_2,e_4,e_2,e_4) \\ 
&- 2 \, R(e_1,e_2,e_3,e_4) = 0 
\end{align*} 
is invariant under parallel transport.
\end{proposition}

\textbf{Proof.} 
Since $(M,g)$ is an Einstein manifold, we have $\text{\rm Ric}_g = \rho \, g$ for some constant $\rho$. Consequently, the metrics $(1 - 2\rho t) \, g$ form a solution to the Ricci flow with nonnegative isotropic curvature. Hence, the assertion follows from Proposition 8 in \cite{Brendle-Schoen2}. \\

\section{K\"ahler-Einstein manifolds}

Let $(M,g)$ be a compact, simply connected Riemannian manifold of dimension $2m \geq 4$ with holonomy group $\text{\rm Hol}(M,g) = U(m)$. Then $(M,g)$ is a K\"ahler manifold. The following theorem was established by S.~Goldberg and S.~Kobayashi:

\begin{theorem}[S.~Goldberg and S.~Kobayashi \cite{Goldberg-Kobayashi}]
\label{goldberg.kobayashi.theorem}
Assume that $(M,g)$ is Einstein. Moreover, suppose that $(M,g)$ has positive orthogonal bisectional curvature; that is, 
\[R(X,JX,Y,JY) > 0\] 
for all points $p \in M$ and all unit vectors $X,Y \in T_p M$ satisfying $g(X,Y) = g(JX,Y) = 0$. Then $(M,g)$ has constant holomorphic sectional curvature.
\end{theorem}

In \cite{Goldberg-Kobayashi}, this result is stated under the stronger assumption that $(M,g)$ has positive holomorphic bisectional curvature (see \cite{Goldberg-Kobayashi}, Theorem 5). However, the proof in \cite{Goldberg-Kobayashi} only uses the condition that $(M,g)$ has positive orthogonal bisectional curvature.

The following result is a consequence of Proposition \ref{borderline} (see also \cite{Seshadri}):

\begin{proposition}
\label{positive.orthogonal.bisectional.curvature}
Assume that $(M,g)$ is Einstein. If $(M,g)$ has nonnegative isotropic curvature, then $(M,g)$ has positive orthogonal bisectional curvature.
\end{proposition}

\textbf{Proof.}
Consider two unit vectors $X,Y \in T_p M$ satisfying $g(X,Y) = g(JX,Y) = 0$. Then 
\begin{align*}
&R(X,Y,X,Y) + R(X,JY,X,JY) \\ 
&+ R(JX,Y,JX,Y) + R(JX,JY,JX,JY) \\ 
&= 2 \, R(X,JX,Y,JY). 
\end{align*} 
Since $(M,g)$ has nonnegative isotropic curvature, it follows that 
\[R(X,JX,Y,JY) \geq 0.\]
It remains to show that $R(X,JX,Y,JY) \neq 0$. To prove this, we argue by contradiction. Suppose that $R(X,JX,Y,JY) = 0$. This implies that the four-frame $\{X,JX,Y,-JY\}$ has zero isotropic curvature. Let us fix a point $q \in M$ and two unit vectors $Z,W \in T_q M$ satisfying $g(Z,W) = g(JZ,W) = 0$. We claim that 
\begin{equation} 
\label{identity.1}
R(Z,JZ,W,JW) = 0. 
\end{equation}
Since $\text{\rm Hol}(M,g) = U(m)$, we can find a piecewise smooth path $\gamma: [0,1] \to M$ such that $\gamma(0) = p$, $\gamma(1) = q$, $P_\gamma X = Z$, and $P_\gamma Y = W$. By Proposition \ref{borderline}, the four-frame $\{P_\gamma X,P_\gamma JX,P_\gamma Y,-P_\gamma JY\}$ has zero isotropic curvature. Consequently, the four-frame $\{Z,JZ,W,-JW\}$ has zero isotropic curvature. Thus, we conclude that $R(Z,JZ,W,JW) = 0$, as claimed.

In the next step, we apply the identity (\ref{identity.1}) to the vectors $\frac{1}{\sqrt{2}} \, (Z + W)$ and $\frac{1}{\sqrt{2}} \, (Z - W)$. This yields 
\begin{align} 
\label{identity.2}
0 &= R(Z + W,JZ + JW,Z - W,JZ - JW) \notag \\ 
&= R(Z,JZ,Z,JZ) + R(W,JW,W,JW) \\
&+ 2 \, R(Z,JZ,W,JW) - 4 \, R(Z,JW,Z,JW). \notag
\end{align}
Similarly, if we apply the identity (\ref{identity.1}) to the vectors $\frac{1}{\sqrt{2}} \, (Z + JW)$ and $\frac{1}{\sqrt{2}} \, (Z - JW)$, then we obtain
\begin{align} 
\label{identity.3}
0 &= R(Z + JW,JZ - W,Z - JW,JZ + W) \notag \\ 
&= R(Z,JZ,Z,JZ) + R(W,JW,W,JW) \\
&+ 2 \, R(Z,JZ,W,JW) - 4 \, R(Z,W,Z,W). \notag
\end{align}
We now take the arithmetic mean of (\ref{identity.2}) and (\ref{identity.3}). This implies 
\begin{equation} 
\label{identity.4}
R(Z,JZ,Z,JZ) + R(W,JW,W,JW) = 0 
\end{equation} 
for all unit vectors $Z,W \in T_q M$ satisfying $g(Z,W) = g(JZ,W) = 0$. 

It follows from (\ref{identity.1}) and (\ref{identity.4}) that the scalar curvature of $(M,g)$ is equal to zero. Since $(M,g)$ has nonnegative isotropic curvature, Proposition 2.5 in \cite{Micallef-Wang} implies that the Weyl tensor of $(M,g)$ vanishes. Consequently, $(M,g)$ is flat. This is a contradiction. \\

Combining Theorem \ref{goldberg.kobayashi.theorem} and Proposition \ref{positive.orthogonal.bisectional.curvature}, we can draw the following conclusion:

\begin{corollary}
\label{Kahler}
Assume that $(M,g)$ is Einstein. If $(M,g)$ has nonnegative isotropic curvature, then $(M,g)$ has constant holomorphic sectional curvature.
\end{corollary}

\section{Quaternionic-K\"ahler manifolds}

Throughout this section, we will assume that $(M,g)$ is a compact, simply connected Riemannian manifold of dimension $4m \geq 8$ with holonomy group $\text{\rm Hol}(M,g) = \text{\rm Sp}(m) \cdot \text{\rm Sp}(1)$. These assumptions imply that $(M,g)$ is a quaternionic-K\"ahler manifold. Hence, there exists a subbundle $\mathcal{G} \subset \text{\rm End}(TM)$ of rank $3$ with the following properties: 
\begin{itemize}
\item $\mathcal{G}$ is invariant under parallel transport.
\item Given any point $p \in M$, we can find linear transformations $I,J,K \in \text{\rm End}(T_p M)$ such that $I^2 = J^2 = K^2 = IJK = -\text{\rm id}$, 
\[g(X,Y) = g(IX,IY) = g(JX,JY) = g(KX,KY)\] 
for all vectors $X,Y \in T_p M$, and 
\[\mathcal{G}_p = \{aI + bJ + cK \in \text{\rm End}(T_p M): a,b,c \in \mathbb{R}\}.\]
\end{itemize}
For each point $p \in M$, we define 
\[\mathcal{J}_p = \{aI + bJ + cK \in \text{\rm End}(T_p M): a,b,c \in \mathbb{R}, \, a^2 + b^2 + c^2 = 1\}.\] 
Note that $\mathcal{J}_p \subset \mathcal{G}_p$ is a sphere of radius $\sqrt{4m}$ centered at the origin. In particular, $\mathcal{J}_p$ is independent of the particular choice of $I,J,K$.

By a theorem of D.~Alekseevskii (see \cite{Besse}, Section 14.41), the curvature tensor of $(M,g)$ can be written in the form $R = R_1 + \kappa \, R_0$ for some constant $\kappa$. Here, $R_1$ is a hyper-K\"ahler curvature tensor; that is, 
\begin{align*} 
R_1(X,Y,Z,W) &= R_1(X,Y,IZ,IW) \\ 
&= R_1(X,Y,JZ,JW) \\ 
&= R_1(X,Y,KZ,KW) 
\end{align*} 
for all vectors $X,Y,Z,W \in T_p M$. Moreover, $R_0$ is defined by 
\begin{align*} 
&4 \, R_0(X,Y,Z,W) \\ &= g(X,Z) \, g(Y,W) - g(X,W) \, g(Y,Z) \\ 
&+ 2 \, g(IX,Y) \, g(IZ,W) + g(IX,Z) \, g(IY,W) - g(IX,W) \, g(IY,Z) \\ 
&+ 2 \, g(JX,Y) \, g(JZ,W) + g(JX,Z) \, g(JY,W) - g(JX,W) \, g(JY,Z) \\ 
&+ 2 \, g(KX,Y) \, g(KZ,W) + g(KX,Z) \, g(KY,W) - g(KX,W) \, g(KY,Z) 
\end{align*}
for all vectors $X,Y,Z,W \in T_p M$. Note that this definition is independent of the particular choice of $I,J,K$.

In the next step, we show that $Q(R) = Q(R_1) + \kappa^2 \, Q(R_0)$. In order to prove this, we need two lemmata:

\begin{lemma}
\label{cross.terms.1}
Fix a point $p \in M$. Let us define an algebraic curvature tensor $S$ on $T_p M$ by 
\[S(X,Y,Z,W) = g(X,Z) \, g(Y,W) - g(X,W) \, g(Y,Z)\] 
for all vectors $X,Y,Z,W \in T_p M$. Then $B(R_1,S) = 0$.
\end{lemma}

\textbf{Proof.}
Let $\{e_1, \hdots, e_{4m}\}$ be an orthonormal basis of $T_p M$. Since the Ricci tensor of $R_1$ vanishes, we have 
\[\sum_{p,q=1}^{4m} R_1(X,Y,e_p,e_q) \, S(Z,W,e_p,e_q) = 2 \, R_1(X,Y,Z,W)\] 
and 
\[\sum_{p,q=1}^{4m} R_1(X,e_p,Z,e_q) \, S(Y,e_p,W,e_q) = -R_1(X,W,Z,Y)\] 
for all vectors $X,Y,Z,W \in T_p M$. Using the first Bianchi identity, we obtain  
\begin{align*} 
B(R_1,S)(X,Y,Z,W) 
&= R_1(X,Y,Z,W) + R_1(Z,W,X,Y) \\ 
&- R_1(X,W,Z,Y) - R_1(Y,Z,W,X) \\ 
&+ R_1(X,Z,W,Y) + R_1(Y,W,Z,X) \\ 
&= 0
\end{align*} 
for all vectors $X,Y,Z,W \in T_p M$. This completes the proof. \\

\begin{lemma}
\label{cross.terms.2}
Fix a point $p \in M$ and an almost complex structure $J \in \mathcal{J}_p$. Let us define an algebraic curvature tensor $S$ on $T_p M$ by 
\begin{align*} 
S(X,Y,Z,W) &= 2 \, g(JX,Y) \, g(JZ,W) \\ 
&+ g(JX,Z) \, g(JY,W) - g(JX,W) \, g(JY,Z) 
\end{align*} 
for all vectors $X,Y,Z,W \in T_p M$. Then $B(R_1,S) = 0$.
\end{lemma}

\textbf{Proof.}
Let $\{e_1, \hdots, e_{4m}\}$ be an orthonormal basis of $T_p M$. Since $R_1$ is a hyper-K\"ahler curvature tensor, we have 
\[\sum_{p,q=1}^{4m} R_1(X,Y,e_p,e_q) \, S(Z,W,e_p,e_q) = 2 \, R_1(X,Y,Z,W)\] 
and 
\begin{align*} 
&\sum_{p,q=1}^{4m} R_1(X,e_p,Z,e_q) \, S(Y,e_p,W,e_q) \\ 
&= 2 \, R_1(X,JY,Z,JW) + R_1(X,JW,Z,JY) 
\end{align*}
for all vectors $X,Y,Z,W \in T_p M$. This implies 
\begin{align*} 
B(R_1,S)(X,Y,Z,W) 
&= R_1(X,Y,Z,W) + R_1(Z,W,X,Y) \\ 
&+ 2 \, R_1(X,JY,Z,JW) + R_1(X,JW,Z,JY) \\ 
&+ 2 \, R_1(Y,JX,W,JZ) + R_1(Y,JZ,W,JX) \\ 
&- 2 \, R_1(X,JY,W,JZ) - R_1(X,JZ,W,JY) \\ 
&- 2 \, R_1(Y,JX,Z,JW) - R_1(Y,JW,Z,JX) 
\end{align*}
for all vectors $X,Y,Z,W \in T_p M$. Using the first Bianchi identity, we obtain 
\begin{align*} 
&B(R_1,S)(X,Y,Z,W) \\ 
&= 2 \, R_1(X,Y,Z,W) + 2 \, R_1(X,JW,Y,JZ) - 2 \, R_1(X,JZ,Y,JW) \\ 
&= 2 \, R_1(X,Y,JZ,JW) + 2 \, R_1(X,JW,Y,JZ) - 2 \, R_1(X,JZ,Y,JW) \\ 
&= 0 
\end{align*} 
for all vectors $X,Y,Z,W \in T_p M$. From this, the assertion follows. \\

\begin{proposition}
\label{cross.terms}
We have $Q(R) = Q(R_1) + \kappa^2 \, Q(R_0)$.
\end{proposition}

\textbf{Proof.} 
Fix a point $p \in M$. Moreover, let $I,J,K \in \mathcal{J}_p$ be three almost complex structures satisfying $IJK = -\text{\rm id}$. We define 
\begin{align*} 
S_0(X,Y,Z,W) &= g(X,Z) \, g(Y,W) - g(X,W) \, g(Y,Z), \\[1.5mm] 
S_1(X,Y,Z,W) &= 2 \, g(IX,Y) \, g(IZ,W) \\ 
&+ g(IX,Z) \, g(IY,W) - g(IX,W) \, g(IY,Z), \\[1.5mm] 
S_2(X,Y,Z,W) &= 2 \, g(JX,Y) \, g(JZ,W) \\ 
&+ g(JX,Z) \, g(JY,W) - g(JX,W) \, g(JY,Z), \\[1.5mm] 
S_3(X,Y,Z,W) &= 2 \, g(KX,Y) \, g(KZ,W) \\ 
&+ g(KX,Z) \, g(KY,W) - g(KX,W) \, g(KY,Z)
\end{align*} 
for all vectors $X,Y,Z,W \in T_p M$. It follows from Lemma \ref{cross.terms.1} and Lemma \ref{cross.terms.2} that 
\[B(R_1,S_0) = B(R_1,S_1) = B(R_1,S_2) = B(R_1,S_3) = 0.\] 
Since $S_0 + S_1 + S_2 + S_3 = 4 \, R_0$, we conclude that $B(R_1,R_0) = 0$. This implies 
\[Q(R) = Q(R_1) + 2\kappa \, B(R_1,R_0) + \kappa^2 \, Q(R_0) = Q(R_1) + \kappa^2 \, Q(R_0),\] 
as claimed. \\

\begin{proposition}
\label{Q.in.Kahler.case}
Fix a point $p \in M$ and an almost complex structure $J \in \mathcal{J}_p$. Moreover, let $\{e_1, \hdots, e_{4m}\}$ be an orthonormal basis of $T_p M$. Then 
\begin{align*} 
Q(R_1)(X,JX,X,JX) 
&\leq -2 \, R_1(X,JX,X,JX)^2 \\ 
&+ 2 \sum_{p,q=1}^{4m} R_1(X,JX,e_p,e_q)^2 
\end{align*}
for every unit vector $X \in T_p M$.
\end{proposition}

\textbf{Proof.} 
By definition of $Q(R_1)$, we have 
\begin{align*} 
Q(R_1)(X,JX,X,JX) 
&= \sum_{p,q=1}^{4m} R_1(X,JX,e_p,e_q)^2 \\ 
&+ 2 \sum_{p,q=1}^{4m} R_1(X,e_p,X,e_q) \, R_1(JX,e_p,JX,e_q) \\ 
&- 2 \sum_{p,q=1}^{4m} R_1(X,e_p,JX,e_q) \, R_1(JX,e_p,X,e_q). 
\end{align*} From this, we deduce that 
\begin{align*} 
Q(R_1)(X,JX,X,JX) 
&= \sum_{p,q=1}^{4m} R_1(X,JX,e_p,e_q)^2 \\ 
&- 4 \sum_{p,q=1}^{4m} R_1(X,e_p,JX,e_q) \, R_1(JX,e_p,X,e_q). 
\end{align*} 
The expression on the right-hand side is independent of the choice of the orthonormal basis $\{e_1, \dots, e_{4m}\}$. Hence, we may assume without loss of generality that $e_1 = X$ and $e_2 = JX$. This implies 
\begin{align*} 
&-4 \sum_{p,q=1}^{4m} R_1(X,e_p,JX,e_q) \, R_1(JX,e_p,X,e_q) \\ 
&= -4 \sum_{p,q=3}^{4m} R_1(X,e_p,JX,e_q) \, R_1(JX,e_p,X,e_q) \\ 
&\leq \sum_{p,q=3}^{4m} (R_1(X,e_p,JX,e_q) - R_1(JX,e_p,X,e_q))^2 \\ 
&= \sum_{p,q=3}^{4m} R_1(X,JX,e_p,e_q)^2 \\ 
&\leq -2 \, R_1(X,JX,X,JX)^2 + \sum_{p,q=1}^{4m} R_1(X,JX,e_p,e_q)^2. 
\end{align*} 
Putting these facts together, the assertion follows. \\

\begin{lemma}
\label{max.hol.sect.curv}
Fix a point $p \in M$ and an almost complex structure $J \in \mathcal{J}_p$. Suppose that $X \in T_p M$ is a unit vector with the property that $R_1(X,JX,X,JX)$ is maximal. Moreover, let $Y \in T_p M$ be a unit vector satisfying $g(X,Y) = g(JX,Y) = 0$. Then 
\[R_1(X,JX,X,Y) = R_1(X,JX,X,JY) = 0\] 
and 
\[2 \, R_1(X,JX,Y,JY) \leq R_1(X,JX,X,JX).\] 
\end{lemma} 

\textbf{Proof.} 
Since $R_1(X,JX,X,JX)$ is maximal, we have 
\[(1 + s^2)^{-2} \, R_1(X + sY,JX + s \, JY,X + sY,JX + s \, JY) \leq R_1(X,JX,X,JX)\] 
for all $s \in \mathbb{R}$. Consequently, we have 
\[\frac{d}{ds} \Big ( (1 + s^2)^{-2} \, R_1(X + sY,JX + s \, JY,X + sY,JX + s \, JY) \Big ) \Big |_{s=0} = 0\]  
and 
\[\frac{d^2}{ds^2} \Big ( (1+s^2)^{-2} \, R_1(X + sY,JX + s \, JY,X + sY,JX + s \, JY) \Big ) \Big |_{s=0} \leq 0.\] 
This implies 
\[R_1(X,JX,X,JY) = 0\] 
and 
\[2 \, R_1(X,JY,X,JY) \leq R_1(X,JX,X,JX) - R_1(X,JX,Y,JY).\] 
Replacing $Y$ by $JY$ yields 
\[R_1(X,JX,X,Y) = 0\]
and
\[2 \, R_1(X,Y,X,Y) \leq R_1(X,JX,X,JX) - R_1(X,JX,Y,JY).\] 
Putting these facts together, we obtain 
\begin{align*} 
R_1(X,JX,Y,JY) 
&= R_1(X,Y,X,Y) + R_1(X,JY,X,JY) \\ 
&\leq R_1(X,JX,X,JX) - R_1(X,JX,Y,JY). 
\end{align*}
From this, the assertion follows. \\

\begin{theorem}
\label{quaternionic.Kahler.1}
Assume that $R_1(X,JX,X,JX) < \kappa$ for every point $p \in M$, every almost complex structure $J \in \mathcal{J}_p$, and every unit vector $X \in T_p M$. Then $R_1$ vanishes identically.
\end{theorem}

\textbf{Proof.}
Note that $R_1$ is a hyper-K\"ahler curvature tensor. Therefore, the Ricci tensor of $R_1$ is equal to $0$. Using the identity $R = R_1 + \kappa \, R_0$, we obtain $\text{\rm Ric}_g = (m+2)\kappa \, g$. Hence, Proposition \ref{curvature.identity} implies that 
\[\Delta R + Q(R) = (2m+4)\kappa \, R.\] 
Since $R_0$ is parallel, we have $\Delta R = \Delta R_1$. Moreover, we have $Q(R_0) = (2m+4) \, R_0$. Using Proposition \ref{cross.terms}, we obtain $Q(R) = Q(R_1) + (2m+4)\kappa^2 \, R_0$. Thus, we conclude that 
\[\Delta R_1 + Q(R_1) = (2m+4)\kappa \, R_1.\] 
By compactness, we can find a point $p \in M$, an almost complex structure $J \in \mathcal{J}_p$, and a unit vector $X \in T_p M$ such that $R_1(X,JX,X,JX)$ is maximal. This implies 
\[(D_{v,v}^2 R_1)(X,JX,X,JX) \leq 0\] 
for all vectors $v \in T_p M$. Taking the trace over $v \in T_p M$ yields 
\[(\Delta R_1)(X,JX,X,JX) \leq 0.\] 
Putting these facts together, we conclude that 
\begin{equation} 
\label{lower.bound.for.Q}
Q(R_1)(X,JX,X,JX) \geq (2m+4)\kappa \, R_1(X,JX,X,JX). 
\end{equation} 

We now analyze the term $Q(R_1)(X,JX,X,JX)$. For abbreviation, let $w_1 = X$ and $w_2 = IX$. We can find vectors $w_3, \hdots, w_{2m} \in T_p M$ such that $\{w_1,Jw_1,w_2,Jw_2, \hdots, w_{2m},Jw_{2m}\}$ is an orthonormal basis of $T_p M$ and 
\[R_1(X,JX,w_\alpha,w_\beta) = R_1(X,JX,w_\alpha,Jw_\beta) = 0\] 
for $3 \leq \alpha < \beta \leq 2m$. It follows from Lemma \ref{max.hol.sect.curv} that 
\[R_1(X,JX,X,w_\beta) = R_1(X,JX,X,Jw_\beta) = 0\] 
for $2 \leq \beta \leq 2m$. Moreover, we have 
\[R_1(X,JX,X,Iw_\beta) = R_1(X,JX,X,JIw_\beta) = 0\] 
for $3 \leq \beta \leq 2m$. This implies 
\[R_1(X,JX,IX,w_\beta) = R_1(X,JX,IX,Jw_\beta) = 0\] 
for $3 \leq \beta \leq 2m$. Putting these facts together, we conclude that 
\begin{equation} 
\label{y.1}
R_1(X,JX,w_\alpha,w_\beta) = R_1(X,JX,w_\alpha,Jw_\beta) = 0 
\end{equation}
for $1 \leq \alpha < \beta \leq 2m$.

Using Lemma \ref{max.hol.sect.curv}, we obtain 
\[2 \, R_1(X,JX,w_\alpha,Jw_\alpha) \leq R_1(X,JX,X,JX)\] 
and 
\[2 \, R_1(X,JX,Iw_\alpha,JIw_\alpha) \leq R_1(X,JX,X,JX)\] 
for $3 \leq \alpha \leq 2m$. The latter inequality implies that 
\[-2 \, R_1(X,JX,w_\alpha,Jw_\alpha) \leq R_1(X,JX,X,JX)\] 
for $3 \leq \alpha \leq 2m$. Thus, we conclude that 
\begin{equation} 
\label{y.2}
4 \, R_1(X,JX,w_\alpha,Jw_\alpha)^2 \leq R_1(X,JX,X,JX)^2 
\end{equation}
for $3 \leq \alpha \leq 2m$. 

By Proposition \ref{Q.in.Kahler.case}, we have 
\begin{align*} 
Q(R_1)(X,JX,X,JX) &\leq -2 \, R_1(X,JX,X,JX)^2 \\ 
&+ 4 \sum_{\alpha,\beta=1}^{2m} R_1(X,JX,w_\alpha,w_\beta)^2 \\ 
&+ 4 \sum_{\alpha,\beta=1}^{2m} R_1(X,JX,w_\alpha,Jw_\beta)^2. 
\end{align*}
Using (\ref{y.1}) and (\ref{y.2}), we obtain 
\begin{align} 
\label{upper.bound.for.Q}
&Q(R_1)(X,JX,X,JX) \notag \\ 
&\leq -2 \, R_1(X,JX,X,JX)^2 + 4 \sum_{\alpha=1}^{2m} R_1(X,JX,w_\alpha,Jw_\alpha)^2 \notag \\ 
&= 6 \, R_1(X,JX,X,JX)^2 + 4 \sum_{\alpha=3}^{2m} R_1(X,JX,w_\alpha,Jw_\alpha)^2 \\ 
&\leq (2m+4) \, R_1(X,JX,X,JX)^2. \notag
\end{align}
Combining (\ref{lower.bound.for.Q}) and (\ref{upper.bound.for.Q}), we conclude that 
\[\kappa \, R_1(X,JX,X,JX) \leq R_1(X,JX,X,JX)^2.\] 
Since $R_1(X,JX,X,JX) < \kappa$, it follows that $R_1(X,JX,X,JX) \leq 0$. Therefore, $R_1$ has nonpositive holomorphic sectional curvature. Since the scalar curvature of $R_1$ is equal to $0$, we conclude that $R_1$ vanishes identically. \\

\begin{proposition}
Assume that $(M,g)$ has nonnegative isotropic curvature. Then $R_1(X,JX,X,JX) < \kappa$ for every point $p \in M$, every almost complex structure $J \in \mathcal{J}_p$, and every unit vector $X \in T_p M$. 
\end{proposition}

\textbf{Proof.}
Fix a point $p \in M$ and a unit vector $X \in T_p M$. Moreover, let $I,J,K \in \mathcal{J}_p$ be three almost complex structures satisfying $IJK = -\text{\rm id}$. For abbreviation, we put $Y = IX$. Then 
\begin{align*}
&R_1(X,Y,X,Y) + R_1(X,JY,X,JY) \\ 
&+ R_1(JX,Y,JX,Y) + R_1(JX,JY,JX,JY) \\ 
&= 2 \, R_1(X,JX,Y,JY). 
\end{align*} 
Moreover, we have 
\begin{align*}
&R_0(X,Y,X,Y) = R_0(X,JY,X,JY) = 1, \\ 
&R_0(JX,Y,JX,Y) = R_0(JX,JY,JX,JY) = 1, \\ 
&R_0(X,JX,Y,JY) = 0 
\end{align*} 
by definition of $R_0$. Using the identity $R = R_1 + \kappa \, R_0$, we obtain 
\begin{align*}
&R(X,Y,X,Y) + R(X,JY,X,JY) \\ 
&+ R(JX,Y,JX,Y) + R(JX,JY,JX,JY) \\ 
&+ 2 \, R(X,JX,Y,JY) \\ 
&= 4 \, (\kappa + R_1(X,JX,Y,JY)) \\ 
&= 4 \, (\kappa - R_1(X,JX,X,JX)). 
\end{align*} 
Since $(M,g)$ has nonnegative isotropic curvature, it follows that 
\[R_1(X,JX,X,JX) \leq \kappa.\] 
It remains to show that $R_1(X,JX,X,JX) \neq \kappa$. To prove this, we argue by contradiction. Suppose that $R_1(X,JX,X,JX) = \kappa$. This implies that the four-frame $\{X,JX,Y,-JY\}$ has zero isotropic curvature. Given any unit vector $Z \in T_p M$, we can find a linear isometry $L: T_p M \to T_p M$ which commutes with $I,J,K$ and satisfies $LX = Z$. Since $\text{\rm Hol}(M,g) = \text{\rm Sp}(m) \cdot \text{\rm Sp}(1)$, there exists a piecewise smooth path $\gamma: [0,1] \to M$ such that $\gamma(0) = \gamma(1) = p$ and $P_\gamma = L$. By Proposition \ref{borderline}, the four-frame $\{P_\gamma X,P_\gamma JX,P_\gamma Y,-P_\gamma JY\}$ has zero isotropic curvature. Hence, if we put $W = IZ$, then the four-frame $\{Z,JZ,W,-JW\}$ has zero isotropic curvature. Consequently, we have 
\[R_1(Z,JZ,Z,JZ) = \kappa\] 
for all unit vectors $Z \in T_p M$. Since $R_1$ is a hyper-K\"ahler curvature tensor, we conclude that $\kappa = 0$. Hence, Proposition 2.5 in \cite{Micallef-Wang} implies that $(M,g)$ is flat. This is a contradiction. \\

\begin{corollary}
\label{quaternionic.Kahler.2}
If $(M,g)$ has nonnegative isotropic curvature, then $R_1$ vanishes identically.
\end{corollary}

\section{Proof of the main theorem}

In this section, we show that every Einstein manifold with nonnegative isotropic curvature is locally symmetric. To that end, we need the following result:

\begin{theorem}
\label{Tachibana.1}
Let $(M,g)$ be a compact Einstein manifold of dimension $n \geq 4$. If $(M,g)$ has positive isotropic curvature, then $(M,g)$ has constant sectional curvature.
\end{theorem}

\textbf{Proof.}
After rescaling the metric if necessary, we may assume that $\text{\rm Ric}_g = (n-1) \, g$. Using Proposition \ref{curvature.identity}, we obtain 
\[\Delta R + Q(R) = 2(n-1) \, R.\] 
We now define 
\[S_{ijkl} = R_{ijkl} - \kappa \, (g_{ik} \, g_{jl} - g_{il} \, g_{jk}),\] 
where $\kappa$ is a positive constant. Note that $S$ is an algebraic curvature tensor. Let $\kappa$ be the largest constant with the property that $S$ has nonnegative isotropic curvature. Then there exists a point $p \in M$ and a four-frame $\{e_1,e_2,e_3,e_4\} \subset T_p M$ such that 
\begin{align*} 
&S(e_1,e_3,e_1,e_3) + S(e_1,e_4,e_1,e_4) \\ 
&+ S(e_2,e_3,e_2,e_3) + S(e_2,e_4,e_2,e_4) \\ 
&- 2 \, S(e_1,e_2,e_3,e_4) = 0. 
\end{align*}
Hence, it follows from Proposition \ref{Q.terms} that 
\begin{align} 
\label{Q.term.1}
&Q(S)(e_1,e_3,e_1,e_3) + Q(S)(e_1,e_4,e_1,e_4) \notag \\ 
&+ Q(S)(e_2,e_3,e_2,e_3) + Q(S)(e_2,e_4,e_2,e_4) \\ 
&- 2 \, Q(S)(e_1,e_2,e_3,e_4) \geq 0. \notag 
\end{align}
We next observe that 
\begin{align*} 
Q(S)_{ijkl} 
&= Q(R)_{ijkl} + 2(n-1) \, \kappa^2 \, (g_{ik} \, g_{jl} - g_{il} \, g_{jk}) \\ 
&- 2\kappa \, (\text{\rm Ric}_{ik} \, g_{jl} - \text{\rm Ric}_{il} \, g_{jk} - \text{\rm Ric}_{jk} \, g_{il} + \text{\rm Ric}_{jl} \, g_{ik}), 
\end{align*} 
hence 
\[Q(S)_{ijkl} = Q(R)_{ijkl} + 2(n-1) \, \kappa \, (\kappa - 2) \, (g_{ik} \, g_{jl} - g_{il} \, g_{jk}).\] 
Substituting this into (\ref{Q.term.1}), we obtain 
\begin{align} 
\label{Q.term.2}
&Q(R)(e_1,e_3,e_1,e_3) + Q(R)(e_1,e_4,e_1,e_4) \notag \\ 
&+ Q(R)(e_2,e_3,e_2,e_3) + Q(R)(e_2,e_4,e_2,e_4) \\ 
&- 2 \, Q(R)(e_1,e_2,e_3,e_4) + 8(n-1) \, \kappa \, (\kappa - 2) \geq 0. \notag
\end{align} 
Since $\{e_1,e_2,e_3,e_4\}$ realizes the minimum isotropic curvature of $(M,g)$, we have
\begin{align*} 
&(D_{v,v}^2 R)(e_1,e_3,e_1,e_3) + (D_{v,v}^2 R)(e_1,e_4,e_1,e_4) \\ 
&+ (D_{v,v}^2 R)(e_2,e_3,e_2,e_3) + (D_{v,v}^2 R)(e_2,e_4,e_2,e_4) \\ 
&- 2 \, (D_{v,v}^2 R)(e_1,e_2,e_3,e_4) \geq 0 
\end{align*} 
for all vectors $v \in T_p M$. Taking the trace over $v \in T_p M$ yields 
\begin{align} 
\label{Laplacian.term}
&(\Delta R)(e_1,e_3,e_1,e_3) + (\Delta R)(e_1,e_4,e_1,e_4) \notag \\ 
&+ (\Delta R)(e_2,e_3,e_2,e_3) + (\Delta R)(e_2,e_4,e_2,e_4) \\ 
&- 2 \, (\Delta R)(e_1,e_2,e_3,e_4) \geq 0. \notag 
\end{align} 
We now add (\ref{Q.term.2}) and (\ref{Laplacian.term}) and divide the result by $2(n-1)$. This implies 
\begin{align*}
&R(e_1,e_3,e_1,e_3) + R(e_1,e_4,e_1,e_4) \\ 
&+ R(e_2,e_3,e_2,e_3) + R(e_2,e_4,e_2,e_4) \\ 
&- 2 \, R(e_1,e_2,e_3,e_4) + 4\kappa \, (\kappa - 2) \geq 0. 
\end{align*} 
On the other hand, we have 
\begin{align*}
&R(e_1,e_3,e_1,e_3) + R(e_1,e_4,e_1,e_4) \\ 
&+ R(e_2,e_3,e_2,e_3) + R(e_2,e_4,e_2,e_4) \\ 
&- 2 \, R(e_1,e_2,e_3,e_4) - 4\kappa = 0. 
\end{align*} 
Since $\kappa$ is positive, it follows that $\kappa \geq 1$. Therefore, $S$ has nonnegative isotropic curvature and nonpositive scalar curvature. By Proposition 2.5 in \cite{Micallef-Wang}, the Weyl tensor of $S$ vanishes. From this, the assertion follows. \\

\begin{proposition}
\label{Tachibana.2}
Let $(M,g)$ be a compact, simply connected Einstein manifold of dimension $n \geq 4$ with $\text{\rm Hol}(M,g) = SO(n)$. If $(M,g)$ has nonnegative isotropic curvature, then $(M,g)$ has constant sectional curvature.
\end{proposition}

\textbf{Proof.}
Suppose that $(M,g)$ does not have constant sectional curvature. By Theorem \ref{Tachibana.1}, there exists a point $p \in M$ and an orthonormal four-frame $\{e_1,e_2,e_3,e_4\} \subset T_p M$ such that 
\begin{align*}
&R(e_1,e_3,e_1,e_3) + R(e_1,e_4,e_1,e_4) \\ 
&+ R(e_2,e_3,e_2,e_3) + R(e_2,e_4,e_2,e_4) \\ 
&- 2 \, R(e_1,e_2,e_3,e_4) = 0. 
\end{align*} 
By assumption, the Weyl tensor of $(M,g)$ does not vanish identically. Hence, we can find a point $q \in M$ and an orthonormal four-frame $\{v_1,v_2,v_3,v_4\} \subset T_q M$ such that $R(v_1,v_2,v_3,v_4) \neq 0$. Since $\text{\rm Hol}(M,g) = SO(n)$, there exists a piecewise smooth path $\gamma: [0,1] \to M$ such that $\gamma(0) = p$, $\gamma(1) = q$, and 
\[v_1 = P_\gamma e_1, \quad v_2 = P_\gamma e_2, \quad v_3 = P_\gamma e_3, \quad v_4 = \pm P_\gamma e_4.\] 
Without loss of generality, we may assume that $v_4 = P_\gamma e_4$. (Otherwise, we replace $v_4$ by $-v_4$.) It follows from Proposition \ref{borderline} that 
\begin{align}
\label{iso.1}
&R(v_1,v_3,v_1,v_3) + R(v_1,v_4,v_1,v_4) \notag \\ 
&+ R(v_2,v_3,v_2,v_3) + R(v_2,v_4,v_2,v_4) \\ 
&- 2 \, R(v_1,v_2,v_3,v_4) = 0. \notag
\end{align} 
Using analogous arguments, we obtain 
\begin{align}
\label{iso.2} 
&R(v_1,v_4,v_1,v_4) + R(v_1,v_2,v_1,v_2) \notag \\ 
&+ R(v_3,v_4,v_3,v_4) + R(v_3,v_2,v_3,v_2) \\ 
&- 2 \, R(v_1,v_3,v_4,v_2) = 0 \notag 
\end{align} 
and 
\begin{align}
\label{iso.3} 
&R(v_1,v_2,v_1,v_2) + R(v_1,v_3,v_1,v_3) \notag \\ 
&+ R(v_4,v_2,v_4,v_2) + R(v_4,v_3,v_4,v_3) \\ 
&- 2 \, R(v_1,v_4,v_2,v_3) = 0. \notag 
\end{align} 
Since $(M,g)$ has nonnegative isotropic curvature, it follows that 
\begin{align*} 
R(v_1,v_2,v_3,v_4) \geq 0, \\ 
R(v_1,v_3,v_4,v_2) \geq 0, \\ 
R(v_1,v_4,v_2,v_3) \geq 0. 
\end{align*} 
Using the first Bianchi identity, we conclude that $R(v_1,v_2,v_3,v_4) = 0$. This is a contradiction. \\

\begin{proposition}
\label{Tachibana.3}
Let $(M,g)$ be a compact, simply connected Einstein manifold of dimension $n \geq 4$ with nonnegative isotropic curvature. Moreover, suppose that $(M,g)$ is irreducible. Then $(M,g)$ is isometric to a symmetric space.
\end{proposition}

\textbf{Proof.} 
Suppose that $(M,g)$ is not isometric to a symmetric space. By Berger's holonomy theorem (see e.g. \cite{Besse}, Corollary 10.92), there are four possibilities:

\textit{Case 1:} $\text{\rm Hol}(M,g) = SO(n)$. In this case, Proposition \ref{Tachibana.2} implies that $(M,g)$ has constant sectional curvature. This contradicts the fact that $(M,g)$ is non-symmetric.

\textit{Case 2:} $n = 2m$ and $\text{\rm Hol}(M,g) = U(m)$. In this case, $(M,g)$ is a K\"ahler manifold. Moreover, by Corollary \ref{Kahler}, $(M,g)$ has constant holomorphic sectional curvature. Consequently, $(M,g)$ is isometric to a symmetric space, contrary to our assumption.

\textit{Case 3:} $n = 4m \geq 8$ and $\text{\rm Hol}(M,g) = \text{\rm Sp}(m) \cdot \text{\rm Sp}(1)$. In this case, $(M,g)$ is a quaternionic-K\"ahler manifold. Moreover, it follows from Corollary \ref{quaternionic.Kahler.2} that $(M,g)$ is symmetric. This is a contradiction.

\textit{Case 4:} $n = 16$ and $\text{\rm Hol}(M,g) = \text{\rm Spin}(9)$. In this case, a theorem of D.~Alekseevskii implies that $(M,g)$ is isometric to a symmetric space (see \cite{Alekseevskii}, Corollary 1, or \cite{Brown-Gray}, Theorem 8.1). Again, this is a contradiction. \\

\begin{theorem}
Let $(M,g)$ be a compact Einstein manifold of dimension $n \geq 4$ with nonnegative isotropic curvature. Then $(M,g)$ is locally symmetric.
\end{theorem}

\textbf{Proof.} 
We first consider the case that $(M,g)$ is Ricci flat. In this case, Proposition 2.5 in \cite{Micallef-Wang} implies that the Weyl tensor of $(M,g)$ vanishes. Consequently, $(M,g)$ is flat.

It remains to consider the case that $(M,g)$ has positive Einstein constant. By a theorem of DeRham (cf. \cite{Besse}, Theorem 10.43), the universal cover of $(M,g)$ is isometric to a product of the form $N_1 \times \hdots \times N_j$, where $N_1, \hdots, N_j$ are compact, simply connected, and irreducible. Since $(M,g)$ is an Einstein manifold, it follows that the factors $N_1, \hdots, N_j$ are Einstein manifolds. Since $(M,g)$ has positive Einstein constant, the manifolds $N_1,\hdots,N_j$ are compact by Myers' theorem. By Proposition \ref{Tachibana.3}, each of the factors $N_1, \hdots, N_j$ is isometric to a symmetric space. Consequently, $(M,g)$ is locally symmetric. \\

We conclude this paper with an analysis of the borderline case in the Micallef-Moore theorem. This result follows from Corollary \ref{quaternionic.Kahler.2} and results established in \cite{Brendle-Schoen2}.

\begin{theorem}
\label{pic.classification}
Let $(M,g_0)$ be a compact, simply connected Riemannian manifold of dimension $n \geq 4$ which is irreducible and has nonnegative isotropic curvature. Then one of the following statements holds:
\begin{itemize}
\item[(i)] $M$ is homeomorphic to $S^n$.
\item[(ii)] $n = 2m$ and $(M,g_0)$ is a K\"ahler manifold.
\item[(iii)] $(M,g_0)$ is isometric to a symmetric space.
\end{itemize}
\end{theorem}

\textbf{Proof.} Suppose that $(M,g_0)$ is not isometric to a symmetric space. Let $g(t)$, $t \in [0,T)$, the unique solution of the Ricci flow with initial metric $g_0$. By continuity, we can find a real number $\delta \in (0,T)$ such that $(M,g(t))$ is irreducible and non-symmetric for all $t \in (0,\delta)$. According to Berger's holonomy theorem (cf. \cite{Besse}, Corollary 10.92), there are four possibilities:

\textit{Case 1:} There exists a real number $\tau \in (0,\delta)$ such that $\text{\rm Hol}(M,g(\tau)) = SO(n)$. In this case, Proposition 8 in \cite{Brendle-Schoen2} implies that $(M,g(\tau))$ has positive isotropic curvature. By a theorem of Micallef and Moore \cite{Micallef-Moore}, $M$ is homeomorphic to $S^n$. 

\textit{Case 2:} $n = 2m$ and $\text{\rm Hol}(M,g(t)) = U(m)$ for all $t \in (0,\delta)$. In this case, $(M,g(t))$ is a K\"ahler manifold for all $t \in (0,\delta)$. Since $g(t) \to g_0$ in $C^\infty$, it follows that $(M,g_0)$ is a K\"ahler manifold. 

\textit{Case 3:} $n = 4m \geq 8$ and $\text{\rm Hol}(M,g(\tau)) = \text{\rm Sp}(m) \cdot \text{\rm Sp}(1)$ for some real number $\tau \in (0,\delta)$. In this case, $(M,g(\tau))$ is a quaternionic-K\"ahler manifold. By Corollary \ref{quaternionic.Kahler.2}, $(M,g(\tau))$ is isometric to a symmetric space. This is a contradiction.

\textit{Case 4:} $n = 16$ and $\text{\rm Hol}(M,g(\tau)) = \text{\rm Spin}(9)$ for some real number $\tau \in (0,\delta)$. By Alekseevskii's theorem, $(M,g(\tau))$ is isometric to a symmetric space (see \cite{Alekseevskii}, Corollary 1, or \cite{Brown-Gray}, Theorem 8.1). This contradicts the fact that $(M,g(\tau))$ is non-symmetric. \\

It is possible to strengthen the conclusion in statement (ii) of Theorem \ref{pic.classification}. To that end, we consider a compact, simply connected K\"ahler manifold which is irreducible and has nonnegative isotropic curvature. By a result of Seshadri, any such manifold is biholomorphic to complex projective space or isometric to a symmetric space (cf. \cite{Seshadri}, Theorem 1.2; see also \cite{Siu-Yau}).

\end{document}